\documentclass[reqno,12pt]{amsart}

\usepackage[margin=2.5cm]{geometry}
\usepackage[T1]{fontenc}
\usepackage[latin1]{inputenc}
\numberwithin{equation}{section}
\usepackage[english]{babel}
\usepackage{amsmath,amssymb,amsthm,amsrefs}
\usepackage{enumitem}
\usepackage{mymacros}
\usepackage{hyperref}
\usepackage{xcolor}
\usepackage{orcidlink}
\usepackage{graphicx}
\usepackage{comment}

\DeclareMathOperator{\capa}{cap}

\DeclareMathOperator{\tr}{tr}
\DeclareMathOperator{\supp}{supp}
\DeclareMathOperator{\hol}{Hol}

\renewcommand{\MR}[1]{}

\title{Carleson measures for Hardy-Sobolev spaces in the Siegel upper half-space}

\DeclareFontFamily{U}{mathx}{}
\DeclareFontShape{U}{mathx}{m}{n}{<-> mathx10}{}
\DeclareSymbolFont{mathx}{U}{mathx}{m}{n}
\DeclareMathAccent{\widehat}{0}{mathx}{"70}
\DeclareMathAccent{\widecheck}{0}{mathx}{"71}

\author[N. Chalmoukis]{N. Chalmoukis \orcidlink{0000-0001-5210-8206}}
\address{Dipartimento di Matematica e Applicazioni, Universit\`a degli
Studi di Milano Bicocca, Via R. Cozzi 55,  20125, Milano, Italy}
\email{nikolaos.chalmoukis@unimib.it}

\author[G. Lamberti]{G. Lamberti \orcidlink{0009-0009-0503-0421}}
\address{Univ. Bordeaux, CNRS, Bordeaux INP, IMB, UMR 5251, F-33400 Talence, France} 
\email{giuseppe.lamberti@math.u-bordeaux.fr}

\thanks{N. Chalmoukis is a member of Indam--Gnampa and partially supported by the Hellenic Foundation for Research and Innovation
(H.F.R.I.) under the ``2nd Call for H.F.R.I. Research Projects to support Faculty Members \&
Researchers'' (Project Number: 4662) and the Indam--Gnampa project CUP\textunderscore E53C22001930001.}


\begin{document}
\begin{abstract}
  We give a capacitary type characterization of Carleson measures for a class of Hardy-Sobolev spaces (also known as weighted Dirichlet spaces) on the Siegel upper half-space, introduced by Arcozzi et al. in \cite{Arcozzi2019}. This answers in part a question raised by the same authors. 
\end{abstract}
\maketitle

\section{Introduction}

It is a basic principle of the so-called ``complex method'' that many properties of a sufficiently regular function defined on $\bR$ correspond to properties of its Poisson extension, and its conjugate harmonic function, on the upper half plane $\bC_+:=\{\zeta =x+iy \in \bC : y>0  \}$. Therefore, certain problems of real analysis can be translated in problems concerning holomorphic functions in $\bC_+$ and vice versa. In the realm of several variables this correspondence takes several forms. A particularly fruitful generalization of this idea in several variables was provided by Stein \cite{Stein1993}. To introduce this point of view let us define the \textit{Siegel domain} 
\[
\cU:= \big\{ \zeta \in \bC^{n+1} : \Im \zeta_{n+1} > \frac{1}{4} \sum_{j=1}^n |\zeta_j|^2 \big\}.
\]
The Siegel domain is a biholomorphic copy  of the unit ball $B_{n+1}$ in $\bC^{n+1}$, via the {\it Caley} map $\cC: B_{n+1} \to \cU$ defined by
\[
\cC(\zeta) = \Big(\frac{2\zeta_1}{1 - \zeta_{n+1}},\dots, \frac{2\zeta_n}{1 - \zeta_{n+1}}, i\frac{1+\zeta_{n+1}}{1-\zeta_{n+1}} \Big), \quad \zeta=(\zeta_1,\zeta_2,\dots,\zeta_{n+1}) \in B_{n+1}.
\]
Thus, the topological boundary of $\cU$,
\[
\partial \cU=\big\{ \zeta \in \bC^{n+1}: \Im \zeta_{n+1} = \frac{1}{4} \sum_{j=1}^n |\zeta_j|^2 \big\},
\]
corresponds, via the Caley map, to the boundary $ \partial B_{n+1}$ of the unit ball with the point $(0,\dots,0,1)$ removed. As we shall shortly see, $\partial \cU$ can be identified in a natural way with the {\it Heisenberg group} $\bH_n$, that is the differentiable manifold $\bC^n \times \bR$ equipped with the group law 
\[
[z,t]\cdot[w,s]=\Big[z+w,t+s- \frac{1}{2}\Im(z \cdot \overline{w}) \Big], \quad [z,t],[w,s] \in \bH_n.
\]
The above structure renders $\bH_n$ a non-abelian, nilpotent Lie group of step two with invariant Haar measure the $(2n+1)$-dimensional Lebesgue measure. The investigation of the interplay between holomorphic functions in $\cU$ and the function theory in $\bH_n$ constitutes part of a program which recently has made significant progress, and of which in the sequel we shall only mention a fraction. 

In \cite{Ogden79} Ogdan and V\'agi proved a Paley-Wiener type representation for the holomorphic Hardy space in a class of domains including the Siegel domain $\cU$. The Hardy space $H^2(\cU)$ is defined as the space of functions $F$, holomorphic in $\cU$, such that 
\begin{equation}\label{defn:Hardy_space} \| F \|_{H^2}^2 : = \sup_{r>0} \int_{\partial \cU} |F(\zeta+r \mathbf{i})|^2 dH_{2n+1}(\zeta) < + \infty,\end{equation}
where  $dH_k$ denotes the real $k$-dimensional Hausdorff measure and $\mathbf{i}=(0,\dots,0,i)$.

More recently,  diverse aspects of spaces of holomorphic functions in the Siegel domain and its generalizations have been studied by Calzi and Peloso, including  Carleson measures \cite{Calzi22} and biholomorphic invariant spaces \cites{Calzi22_b, Calzi2023B} (see also \cites{Calzi21, Calzi2023}). In \cite{Arcozzi2019} the authors introduced a class of Hardy-Sobolev spaces in the Siegel domain and proved Paley-Wiener type representation theorems, while in \cite{Arcozzi2021_b} a von Neumann type inequality for unbounded tuples of operators was proved using the Paley-Wiener theorems developed in \cite{Arcozzi2019}.
In this work we are interested in a class of holomorphic Hilbert spaces introduced in \cite{Arcozzi2019}*{Definition p. 1961}. We recall here the basic definitions. Let  $\rho:\bC^{n+1} \to \bR$  be the defining function of the Siegel domain;
\[
\rho(\zeta):= \Im(\zeta_{n+1}) - \frac{1}{4}|\zeta'|^2, \quad \zeta=(\zeta',\zeta_{n+1})\in \bC^{n}\times \bC  . 
\]
We denote by $\hol(\cU)$ the Fr\'echet space of holomorphic function in the Siegel domain and $\hol(\overline{\cU})$ the space of holomorphic functions in an open neighborhood of $\cU$. Following \cite{Arcozzi2019}, we shall say that a function $F\in \hol(\cU)$ {\it vanishes at infinity} if for all $R>0$
\[ \lim_{\Im (\zeta_{n+1}) \to \infty} \sup_{|\zeta'| \leq R} |F(\zeta',\zeta_{n+1})| = 0.    \]

\begin{defn}[Hardy-Sobolev spaces]
Let $0 \leq \alpha < \frac{n+1}{2}$ and let $m$ be an integer such that $m>\alpha$. We define \footnote{We choose to work with a different parametrization of the spaces $H^2_\alpha$ with respect to the one introduced in the original work of Arcozzi et al. The relation between the two is $\alpha = - \frac{\nu+1}{2}$.} the {\it  Hardy-Sobolev } space $H^2_{\alpha}$ as the space of functions $F\in \hol(\cU)$ which vanish at infinity and 
\begin{equation}\label{defn:Hardy_Sobolev}
\int_{\cU} |\rho^m(\zeta)\partial^m_{\zeta_{n+1}} F(\zeta)|^2 \rho^{-(2\alpha+1)}(\zeta) dH_{2n+2}(\zeta)< + \infty.
\end{equation}
 
\end{defn}
It turns out that a function $F\in \hol(\cU)$ satisfies \eqref{defn:Hardy_Sobolev} for some natural $m>\alpha$ if and only if it satisfies it for every natural $m>\alpha$, which justifies the notation $H^2_\alpha$ for the space of such functions. In fact, much more is true. In \cite{Arcozzi2019}*{Theorem 2} the authors proved that for functions in the Hardy-Sobolev space $H^2_{\alpha}$ the quantity 
\begin{equation}
    \| F \|_{H^2_\alpha}^2 := 4^{-m} \Gamma(2m-2\alpha)  \int_{\cU} |\rho^m(\zeta)\partial^m_{\zeta_{n+1}} F(\zeta)|^2 \rho^{-(2\alpha+1)}(\zeta) dH_{2n+2}(\zeta),
\end{equation}
 does not depend on $m\in \bN$ as long as $m>\alpha$ and it is a norm on $H^2_\alpha$. This norm renders the space a Hilbert space with inner product which we denote by $\langle \cdot, \cdot \rangle_{H^2_\alpha}$. In fact, for $\alpha = 0$, this is exactly the Hardy space $H^2(\cU)$ defined in \eqref{defn:Hardy_space}. 
  Furthermore, the spaces $H^2_{\alpha}$ are {\it reproducing kernel Hilbert spaces} with reproducing kernel  given by

\[ K_{\alpha}(\omega, \zeta)  = K_{\alpha}^\zeta(\omega):=\frac{\Gamma(n+1-2\alpha)}{(4\pi)^{n+1}{\Big( \dfrac{\omega_{n+1}-\overline{\zeta_{n+1}}}{2i} - \frac{1}{4} \omega' \cdot \overline{\zeta'} \Big)^{  n + 1-2\alpha } }}, \, (\zeta',\zeta_{n+1}),(\omega',\omega_{n+1}) \in \bC^{n}\times \bC. \]

In the same work Arcozzi et al. raised the question of characterizing {\it Carleson measures } for the class of spaces $H^2_\alpha$. Let us first recall the definition of Carleson measures in our setting.

\begin{defn}[Carleson measure] A positive Radon measure $ \mu $ on  $\cU$ is a Carleson measure for $H^2_\alpha$ if there exists a constant $C>0$, depending only on $\mu$, such that 
\begin{equation}\label{defn:Carleson_measure}
\int_{\cU}|f(\zeta)|^2 d\mu(\zeta) \leq C \norm{f}^2_{H^2_\alpha}, \quad \forall f \in H^2_\alpha.
\end{equation}
\end{defn}
For the family of weighted Bergman spaces in $\cU$ \cite[Definition p. 1961]{Arcozzi2019}, a characterization of their Carleson measures was obtained, in a much more general setting, in \cite{Calzi22}. In \cite{Hormander67} H\"ormander characterized Carleson measures for the Hardy space defined on bounded pseudoconvex domains. We should also mention the characterization of Carleson measures for Bergman spaces in the unit ball obtained in \cites{Cima82, Duren07} and the more recent works \cites{Abate11, Abate19} which extends the previous results to some classes of pseudoconvex domains.

The study of such measures is relevant in many areas of complex and harmonic analysis. For example, Carleson measures have a central role in the characterization of universal interpolating sequences for reproducing kernel Hilbert spaces of holomorphic functions (see \cite{Agler2002} and \cite{Seip04} for more details). Furthermore, they are linked to the multipliers of reproducing kernel Hilbert spaces \cite{Stegenga80}.

In this work we are interested in finding a characterization of Carleson measures for the spaces $H^2_\alpha$ when $ \alpha \in (\frac{n}{2}, \frac{n+1}{2})$. Before discussing our results in more depth, we should mention that analogous results have been obtained first by Stegenga (\cite{Stegenga80}) for the Dirichlet space in the unit disc and then by Ahern and Cohn (\cite{Ahern89}) in the unit ball of $\bC^n$. In the aforementioned work the authors characterized exceptional sets and Carleson measures for Hardy-Sobolev spaces defined in the unit ball of $\bC^n$, i.e. spaces of holomorphic functions $f$ such that the fractional derivative of $f$, up to a certain order, belongs to the Hardy space of the ball. With this aim they developed a potential theory on $\partial B_{n+1}$. In our setting, in a similar manner, we exploit the Riesz potential theory in the Heisenberg group, in order to characterize Carleson measures for the spaces $H^2_\alpha$. It is worth pointing out that since Riesz potential theory on the Heisenberg group is quite well understood, our approach exhibits less technical difficulties and it offers the advantage of providing a link between the Hardy-Sobolev spaces $ H^2_\alpha $ and the potential theory associated to the Kohn Laplacian of the Heisenberg group.

Let us now discuss in some detail the objects that we have introduced. We start by defining an action of $\bH_n$ on the boundary of the Siegel domain. To an element $[z,t] \in \bH_n$, we associate the affine self-mapping of $\overline{\cU}$
\[
L_{[z,t]}:(\zeta',\zeta_{n+1}) \mapsto \big( \zeta'+z,\zeta_{n+1}+t+\frac{i}{2} \zeta'\cdot \overline{z}+\frac{i}{4}|z|^2 \big), \quad (\zeta',\zeta_{n+1}) \in \overline{\cU}.
\]
As it can be readily verified, this is a faithful and simply transitive Lie group action on $\partial \cU$. Furthermore, for every $\zeta \in \overline{\cU}$ and $[z,t] \in \bH_n$ we have that $\rho(\zeta) = \rho(L_{[z,t]}(\zeta))$. This allows us to identify the Heisenberg group with $\partial\cU$ via its action on the origin.

We will use the following parametrization of $\cU$ by means of a foliation of copies of the boundary. We set $\mathbf{U}: =\bH_n \times(0,+\infty)$. Given $\zeta=(\zeta',\zeta_{n+1}) \in \cU$, we define $\Psi(\zeta',\zeta_{n+1}):=[z,t,h] \in \mathbf{U}$ by
\begin{equation}\label{defn:coordinates}
\Bigg\{
\begin{array}{ll}
     z & = \zeta'  \\
     t & = \Re(\zeta_{n+1}) \\
     h & = \Im(\zeta_{n+1}) - \frac{1}{4}|\zeta'|^2.
\end{array}
\end{equation}
Then $\Psi: \overline{\cU} \to \overline{\mathbf{U}}$ is a $C^\infty$-diffeomorphism and $\Psi^{-1}$ is given by
\begin{equation} \label{defn:inverse_coordinates}
\Psi^{-1}[z,t,h] = \big(z, t + \frac{i}{4}|z|^2+ih) =: (\zeta',\zeta_{n+1}).
\end{equation}
We consider the following function
\[
d([z,t,h])=\big(\frac{1}{16}(|z|^2+h)^2+t^2\big)^{\frac{1}{4}}, \quad [z,t,h] \in \overline{\mathbf{U}},
\]
and we let 
\[
d([z,t,h],[w,s,k]) := d([[z,t]\cdot[w,s]^{-1},4(h+k)]).
\]
Restricted on $\bH_n$ this is exactly the Folland-Kaplan gauge of the Heisenberg group. The associated distance function is given by
\[
d([z,t],[w,s]):= d([z,t]\cdot[w,s]^{-1}), \quad [z,t],[w,s]\in \bH_n.
\]

It is worth mentioning that the topology induced by the metric $d$ on $\bH_n$ is equivalent to the Euclidean metric of $\bC^n \times \bR$. We write $B([z,t],r)$ for the open ball centered at $[z,t] \in \bH_n$ and radius $r>0$, with respect to the metric $d$. We have that, for some constant $c_n>0$,
\begin{equation} \label{measure:hball}
H_{2n+1} (B([z,t],r)) = c_n r^{2n+2},
\end{equation}
as can be found for example in \cite{Stein1993}*{Chapter XII, Section 2.5.2}.

Finally , if $F$ is a function defined on $\overline{\cU}$ and $f$ a function defined on $\overline{\mathbf{U}}$, we set $\widehat{F} := F \circ \Psi^{-1}$  and $ \widecheck{f} :=  f \circ \Psi$. Notice that, since $|\det J \Psi|=1,$ we have

\[
\int_{\cU} F(\zeta)dH_{2n+2}(\zeta)= \int_0^\infty \int_{\bH_n} \widehat{F}[z,t,h] dH_{2n+1}[z,t]dh.
\]

In the next definition we denote by $\capa_\alpha$ the $\alpha$ - {\it   Riesz capacity } of a Borel subset $\bH_n$. A discussion of the potential theory on the Heisenberg group will be given in Section \ref{section:potential}

\begin{defn}[Subcapacitary measures]
Let $E$ be a Borel subset of  $\bH_n$, we define the tent based on $E$ as
\[
T(E):=\Psi^{-1} \big\{ [z,t,h] \in \cU : B([z,t],h^{1/2}) \subseteq E \big\}.
\]
Let now $\mu$ be a positive Radon measure on $\cU$. We say that $\mu$ is $\alpha$-\textit{subcapacitary}, $\alpha > 0 $, if there exists a constant $C>0$ such that, for all disjoint collection of balls $B(x_i,r_i), i = 1,\dots, d$ in $\bH_n$, the following inequality holds true  
\begin{equation}\label{defn:subcapacitary}
\sum_{i=1}^d\mu(T(B(x_i,r_i))) \leq C \capa_{\alpha}( \bigcup_{i=1}^d B\big(x_i,r_i)\big).
\end{equation}
\end{defn}

It should be mentioned that the regularity properties of $\capa_\alpha$ (see \cite{Vodopyanov88}), immediately imply that the subcapacitary condition implies that $\mu(T(E)) \leq C \capa_\alpha(E)$ for all Borel sets $E \subseteq \bH_n$.

We can now formulate our main result.
\begin{thm} \label{main_thm}
Let $\mu$ be a positive Radon measure on $\cU$ and $\frac {n}{2} < \alpha < \frac{n+1}{2}$. Then $\mu$ is a Carleson measure for $H^2_\alpha$, if and only if it is $\alpha$-subcapacitary.
\end{thm}

This provides a fairly concrete characterization of Carleson measures, compatible with the one of Ahern and Cohn (\cite{Ahern89}) for the unit ball.

\section{Preliminaries}

\subsection{Elements of Fourier analysis on the Heisenberg group}

In this section we shall introduce the absolutely necessary elements of Fourier analysis in the Heisenberg group. It is a tool that, although not indispensable for our purposes, will simplify some calculations later. For a more detailed exposition we refer the reader to \cite{thangavelu_2023}.

Let $\lambda >0 $ and consider the {\it Fock space} $\cF^\lambda$
\[ \cF^\lambda : = \big\{ F\in \hol(\bC^n) : \Big( \frac{|\lambda|}{2\pi} \Big)^n \int_{\bC^n} |F(z)|^2 e^{-\frac{\lambda}{2}|z|^2}dH_{2n}(z) < + \infty \big\},  \]
while for $\lambda<0$ set $\cF^\lambda : = \cF^{-\lambda}$. The {\it Bargman representation} of the Heisenberg group is then defined as follows. If $[z,t]\in \bH_n$ and $\lambda>0$ we define the unitary operator  $\sigma_\lambda[z,t]$ on $\cF^\lambda$ as a weighted shift operator 
\[ \sigma_\lambda[z,t]F(w): = e^{i\lambda t - \frac{\lambda}{2}w\cdot \overline{z}-\frac{\lambda}{4}|z|^2}F(z+w), \quad w \in \bC^n \]
and $\sigma_{-\lambda}[z,t]:=\sigma_\lambda[\overline{z},-t]$. For a function $f\in L^1(\bH_n)$, we define the {\it Fourier transform} of $f$ as the operator $\sigma_\lambda(f)$ on $\cF^\lambda$, defined by
\[ \sigma_\lambda(f)F(w) = \int_{\bH_n} f[z,t] \sigma_\lambda[z,t] F(w) dH_{2n+1}[z,t].  \]

We denote by $\tr$ the trace norm of a trace class operator on a Hilbert space.
\begin{thm}[Plancherel's Theorem for $\bH_n$]
    Let $f,g \in L^2(\bH_n)$, then the operators $\sigma_\lambda(f),\sigma_\lambda(g)$ are Hilbert-Schmidt operators for a.e. $\lambda \in \bR $ and 
    \begin{equation}
        \int_{\bH_n}f[z,t]\overline{g[z,t]}dH_{2n+1}[z,t] = \frac{1}{(2\pi)^{n+1}} \int_{-\infty}^\infty \tr (\sigma_\lambda(f) \sigma_\lambda(g)^*) |\lambda|^n d\lambda.
    \end{equation}
\end{thm}

Finally we need the following formula for the Fourier transform of the kernel vectors, obtained in \cite{Arcozzi2019}. Write $\zeta, \omega \in \overline{\cU}$ using the $\overline{\mathbf{U}}$-coordinates $[z,t,h], [w,s,k]$, and let $0<\alpha<\frac{n+1}{2}$. Consider the function  $f_{[z,t,h]}[w,s]:=\widehat{K_\alpha}([w,s,0],[z,t,h])$.  It holds that \cite[p. 1978]{Arcozzi2019} for every $\lambda > 0, \sigma_\lambda(f_{[z,t,h]}) = 0$ and for $\lambda <  0$ 
\begin{equation}\label{eq:Fourier_kernel} \sigma_\lambda(f_{[z,t,h]}) =  2^{-2\alpha} e^{h \lambda} |\lambda|^{-2\alpha} P_0\sigma_\lambda[z,t], \end{equation}
where $P_0:\cF^\lambda \to \cF^\lambda$  is the orthogonal projection on the constant functions. Furthermore, one can show (see \cite[p. 1978]{Arcozzi2019}) that
\begin{equation}\label{eq:tr_Heisenberg}
    \tr(P_0\sigma_\lambda[z,t] P_0 \sigma_\lambda [w,s]^*) = e^{-|\lambda|(\frac{1}{4}|w-z|^2-i(s-t+\frac{1}{2}\Im (w \cdot \overline{z}))}, \quad [z,t], [w,s] \in \bH_n.
\end{equation}

\subsection{Potential theory.} \label{section:potential}
In this section we shall briefly recall some aspects of the potential theory in the Heisenberg group associated to Riesz potentials. One can consult \cite{Adams1996} for a general potential theory of which the one we present is an instance.

For $\alpha > 0 $, we define the {\it Riesz kernel} of order $\alpha$ as the function 
\[ I_\alpha(x,y) := \frac{\Gamma(n+1-\alpha)}{2^\alpha (2\pi)^{n+1}d(x,y)^{2n+2-2\alpha}}, \quad x \in \bH_n. \]Notice that, for $\alpha =1 $, this is a scalar multiple of the fundamental solution of the subelliptic Laplacian in the Heisenberg group (see \cite{Folland1973}). 
Given a positive Radon measure on $\bH_n$, we define its $\alpha$-Riesz potential to be the  convolution with the Riesz kernel $I_\alpha$, that is
\[ I_\alpha (\mu) (x) : = \int_{\bH_n} I_\alpha(x,y)d\mu(y), \]
which can be considered as an everywhere defined function taking values in $[0,\infty]$. 
Let now $A\subseteq \bH^n$, we define its (outer) capacity to be the quantity 
\[ \capa_\alpha(A):= \inf\{ \| f \|_{L^2(\bH_n)}^2 : f\geq 0, \,\, f\in L^2(\bH_n), \,\, I_\alpha (f) \geq 1 \,\, \text{on}\,\, A  \}. \]
Restricted to Borel sets, $\capa_\alpha$ is actually inner and outer regular \cite{Adams1996}*{Proposition 2.3.12}. A particularly useful application of the celebrated minimax theorem of von Neumann, allows us to express the capacity of compact sets in a ``dual'' sense \cite{Adams1996}*{Corollary 2.5.2}.

\begin{thm} Let $A\subseteq \bH_n$ be a compact set and $\cM_+(A)$ the collection of positive finite Borel measures on $A$. Then, 
\[
\capa_{\alpha}(A) =\sup \{ \mu(A): \mu \in \cM_+(A), I_\alpha(I_\alpha(\mu))(x) \leq 1 \,\, \text{on} \,\, \supp \mu \}.
\]\end{thm}

We recall the following Theorem (\cite{Adams1996}*{Theorem 2.3.10, 2.5.3}).

\begin{thm} \label{AdamsThm2.5.3}
Let $A \subseteq \bH_n$ be compact. Then there is a finite positive Borel measure $\mu^A$ on $A$, such that $I_\alpha(I_\alpha(\mu^A))(x) \leq 1$ for all $x\in \supp \mu$, $I_\alpha(I_\alpha(\mu^A))(x) \geq 1 $ for $H_{2n+1}$- a.e. $x \in A$. Furthermore,
\[
\mu^A(A)= \int_{\bH_n} (I_\alpha (I_\alpha( \mu^A))^2 dH_{2n+1} = \capa_{\alpha}(A).
\]
\end{thm}

The results concerning potential theory that we have discussed up to this point apply in very general situations. The next theorem of Vodopyanov \cite{Vodopyanov88}*{Theorem 2} rests on more intrinsic properties of the Riesz potentials. 

\begin{thm}[Strong capacitary inequality] \label{strong_capacitary_inequality}
   There exists a positive constant $C$ depending only on $\alpha$ and $n$ such that for every positive function $ f \in L^2(\bH_n)$ we have 
\[
    \int_0^\infty \capa_{\alpha}\big( \{ [z,t] \in \bH_n: I_\alpha (f) [z,t] > \lambda \} \big) \lambda \, d\lambda \leq C \, \|f\|_{L^2(\bH_n)}^2.
\]
\end{thm}

In our setting, the connection between Riesz potentials in the Heisenberg group and the Hardy-Sobolev spaces $H^2_\alpha$ comes from the fact that, as a straightforward computation shows, 
\begin{equation}\label{eq:kernel_riesz}
|\widehat{K_{\frac{\alpha}{2}}}([z,t],[w,s])|= I_\alpha([z,t], [w,s]), \quad [z,t],[w,s] \in \bH_n.
\end{equation}

\subsection{Other preliminaries}

Finally we recall some basic objects and notions from multivariable complex theory. The \textit{Poisson kernel} of the Siegel domain is the positive kernel defined as follows 
\[
P(\zeta, \omega):= \frac{\big( \Im (\zeta_{n+1}) - \frac{1}{4} |\zeta'|^2 \big)^{n+1}}{\left| \frac{\zeta_{n+1} - \overline{\omega_{n+1}}}{2i} - \frac{1}{4} \zeta' \cdot \overline{\omega'} \right|^{2(n+1)}}, \quad \zeta \in \cU, \omega \in \partial \cU.
\]
For a function $f$ defined on $\partial \cU $ such that 
\begin{equation}\label{condition_poisson}
    \int_{\partial \cU } \frac{|f(\zeta)|}{1+\widecheck{d}(\zeta)^{4(n+1)}} dH_{2n+1}(\zeta) < + \infty,
\end{equation}
we can define its {\it Poisson extension} by
\[
P[f](\zeta) := \int_{\partial \cU} P(\zeta,\omega) f(\omega) dH_{2n+1}(\omega).
\]

Now we recall the definition of a admissible region and the admissible maximal function. This will be used later on the proof of the main theorem. Let $\omega = \Psi^{-1}[w,s] \in \partial\mathcal{U}$ and $ \gamma > 1$. The \textit{admissible region} centered at $\omega$ of aperture $\gamma >0$ is defined as
\begin{align*}
\Gamma_\gamma(\omega)  : & = \big\{ \zeta \in \mathcal{U} : \left| \frac{\zeta_{n+1} - \overline{\omega_{n+1}}}{2i} - \frac{1}{4} \zeta' \cdot \overline{\omega'} \right| < \gamma \big( \Im (\zeta_{n+1}) - \frac{1}{4} |\zeta'|^2 \big) \big\} \\ 
& = \Psi^{-1} \big\{ [z,t,h] \in \mathbf{U} : d([z,t,h],[w,s]) < \sqrt{2  \gamma h} \big\}.
\end{align*}

In a natural way one defines also the admissible maximal function. Let $F$ be an holomorphic function in $\mathcal{U}$. The \textit{admissible maximal function} is
\[
M_\gamma F(\omega):= \sup \big\{ |F(\zeta)| : \zeta \in \Gamma_\gamma(\omega) \big\}.
\]

Finally, let $f$ be locally integrable function $\bH_n$, we define the Hardy-Littlewood maximal function as
\[
Mf ([z,t]):= \sup_{r > 0} \frac{1}{r^{2n+2}} \int_{B([z,t],r)} |f[w,s]| dH_{2n+1}[w,s].
\]

The next Lemma is a well known result about maximal functions, see for instance \cite{Stein1993}*{Chapter XIII, Section 7.11}.

\begin{lem} \label{maximal}
For every $\gamma > 1$ there exists $C> 0$, depending only on $\gamma$, such that
\[
M_\gamma P [\widecheck{f} \, ] \leq C \widecheck{M f},
\]
for every positive measurable function $f$ on $\bH_n$.
\end{lem}


\subsection{Notation}
 With the letter $C$ we will denote a constant which might change from appearance to appearance and it depends only on the parameters $n, \alpha, \gamma$. If $f, g$ are positive expressions we write $f \lesssim g $ if there exists a constant $C>0$ as before, such that $f \leq C g$.  Furthermore we will write $f \simeq g$ instead of $f \lesssim g $ and $ g \lesssim f$.

\section{Proof of the main theorem}

We will first prove that the subcapacitary condition is sufficient in order to have the Carleson embedding. We have divided the proof in a sequence of lemmas. The basic idea is that using the  fractional differentiation operator that we construct in Lemma \ref{lem:fractional_diff}, we can ``embed'' the Hardy-Sobolev space, via the Poisson extension, into the space of Riesz potentials on the Heisenberg group (Lemma \ref{AhernLem1.7}). The remaining steps, involving the admissible maximal function and Vodopyanov's strong capacitary inequality, are quite standard.

\begin{lem}\label{lem:fractional_diff}
    Let $0<\alpha<\frac{n+1}{2}$. Then there exists an isometry  $\mathcal{R}^\alpha : H^2_{\alpha} \to H^2$ such that for every $\zeta \in \cU$ we have 
    \begin{equation}\label{prop:fractional_diff}
        \mathcal{R}^\alpha K^\zeta_{\alpha} = 2^{\alpha}K^\zeta_{\frac{\alpha}{2}}.
    \end{equation}
\end{lem}

\begin{proof}
    We claim that it is sufficient to prove that for every  $\alpha \in (0, \frac{n+1}{2}) $, it holds that $K^\zeta_{\frac{\alpha}{2}} \in H^2(\cU) $, and also
    \begin{equation}
        \langle K^{\zeta}_{\frac{\alpha}{2}}, K^{\omega}_{\frac{\alpha}{2}} \rangle_{H^2} =2^{-2\alpha} K_{\alpha}(\omega,\zeta).
    \end{equation}
If this is the case we initially define $\cR^\alpha$ on finite linear combinations of kernel vectors via \eqref{prop:fractional_diff}. Then for $c_1,\dots c_J \in \bC$, $\zeta_1, \dots \zeta_J \in \cU$ we have 
\begin{align*}
\big\| \cR^\alpha(\sum_{j=1}^J c_j K^{\zeta_j}_{\alpha}) \big\|_{H^2}^2 & = 2^{2\alpha} \sum_{j,k=1}^J c_j \overline{c_k}   \langle K^{\zeta_j}_{\frac{\alpha}{2}}, K^{\zeta_k}_{\frac{\alpha}{2}} \rangle_{H^2} \\ 
& =  \sum_{j,k=1}^J  c_j \overline{c_k} K_\alpha(\zeta_k,\zeta_j)  \\
& =  \big\| \sum_{j=1}^J c_j K^{\zeta_j}_{\alpha} \big\|_{H^2_\alpha}^2. 
\end{align*} 
Hence $\cR^\alpha$ extends in a unique way to a Hilbert space isometry. 

Let now prove the initial claim. For $\zeta, \omega \in \cU, $ with $\mathbf{U}$-coordinates $[z,t,h],[w,s,k]$ and $C>0$, a constant depending only on $n, \alpha$, we have 
\begin{align*}
    \| K_{\frac \alpha2}^\zeta \|_{H^2}^2&  =  C \sup_{h>0} \int_{\bH_n} \frac{ dH_{2n+1}[w,s]}{\big( (k+h+\frac{1}{4}|w-z|^2)^2 +(s-t+ \frac{1}{2}\Im(w\cdot \overline{z})^2 \big)^{n+1-\alpha}} \\
    & = C \int_{\bH_n}  \frac{ dH_{2n+1}[w,s]}{\big( (k+\frac{1}{4}|w-z|^2)^2 +(s-t+\frac{1}{2}\Im(w\cdot \overline{z})^2 \big)^{n+1-\alpha}} \\
    &= C \int_{\bH_n} \frac{dH_{2n+1}[w,s]}{\big( (k+\frac{1}{4}|w|^2)^2 +s^2 \big)^{n+1-\alpha}} < + \infty.
\end{align*}
The last quantity is finite for every $k > 0 $ because $\alpha<\frac{n+1}{2}$. 

Next, let $ f_{[z,t,h]} $ as defined in \eqref{eq:Fourier_kernel}, and apply successively Plancherel's theorem and equations \eqref{eq:Fourier_kernel} and \eqref{eq:tr_Heisenberg} to obtain  
\begin{align*}
    \langle K^\zeta_{\frac{\alpha}{2}},K^\omega_{\frac{\alpha}{2}} \rangle_{H^2}  & =\frac{1}{(2\pi)^{n+1}} \int_{-\infty}^0  \tr\big( \sigma_\lambda(f_{[z,t,h]}) \sigma_\lambda(f_{[w,s,k]})^*\big) |\lambda|^n d\lambda \\
& = \frac{2^{-2\alpha}}{(2\pi)^{n+1}} \int_{-\infty}^0  e^{-|\lambda|(h+k)}\tr\Big( P_0\sigma_\lambda[z,t] P_0\sigma_\lambda[w,s]^* \Big) |\lambda|^{n-2\alpha} d\lambda \\
& = \frac{2^{-2\alpha}}{(2\pi)^{n+1}} \int_0^\infty  e^{-\lambda\big(h+k+\frac{1}{4}|w-z|^2-i\big(s-t+\frac{1}{2}\Im(w\cdot\overline{z}) \big)\big)} \lambda^{n-2\alpha} d\lambda \\ 
& = \frac{\Gamma(n+1-2\alpha)}{2^{2\alpha}(2\pi)^{n+1}} \big(h+k+\frac{1}{4}|w-z|^2-i\big(s-t+\frac{1}{2}\Im(w\cdot\overline{z}) \big)^{2\alpha-(n+1)} \\
& = 2^{-2\alpha} K_\alpha(\omega,\zeta).
\end{align*}

This concludes the proof. 
\end{proof}

\begin{lem} \label{AhernLem1.7}
Let $0<\alpha< \frac{n+1}{2}$. Then, for every $F \in H^2_\alpha \cap \hol(\overline{\cU})$, there exists $f \in L^2(\bH_n), f\geq 0$, satisfying $\|F\|_{H^2_\alpha} = \|f\|_{L^2(\bH_n)}$  and also 
\[
|F(\zeta)| \leq 2^\alpha P[\widecheck{I_\alpha (f)}] (\zeta), \quad \forall \zeta \in \cU.
\]
\end{lem}

\begin{proof}
Let $\zeta \in \cU, \Psi(\zeta)=[z,t,h]$. We use the properties of the operator $\cR^\alpha$ defined in Lemma \ref{lem:fractional_diff} and the fact that functions in $H^2(\cU)$ have boundary values almost everywhere \cite{Ogden79}, to obtain
\begin{align} \label{eq:F_formula}
    F(\zeta) = \langle F, K_\alpha^\zeta \rangle_{H^2_\alpha} = \langle \mathcal{R}^\alpha F, \mathcal{R}^\alpha K_\alpha^\zeta \rangle_{H^2} 
    = 2^\alpha \int_{\partial\cU} \mathcal{R}^\alpha F (\omega) K_{\frac{\alpha}{2}}(\zeta,\omega) dH_{2n+1}(\omega).
\end{align}

Since $K^\omega_{\frac{\alpha}{2}}$ is a pluriharmonic function, then we can write it as the Poisson kernel of its boundary values (see \cite{Rudin80}). This gives 
\[ K_\frac{\alpha}{2}(\zeta,\omega) = \int_{\partial \cU} P(\zeta,\eta) K_{\frac{\alpha}{2}}(\eta,\omega) dH_{2n+1}(\eta). \]
Hence substituting this Poisson representation in the expression \eqref{eq:F_formula} and setting $f:= |\widehat{\cR^\alpha F}  |$  we obtain 
\begin{align*}|F(\zeta)| & \leq 2^\alpha \int_{\partial \cU} P(\zeta,\eta) \int_{\partial \cU} | \cR^\alpha F (\omega) \, K_{\frac{\alpha}{2}}(\eta,\omega)| dH_{2n+1}(\omega) dH_{2n+1}(\eta) \\ 
& =  2^\alpha \int_{\partial \cU} P(\zeta,\eta) \int_{\bH_n} |\widehat{\cR^\alpha F} [w,s] | | \widehat{K_{\frac{\alpha}{2}}}(\Psi(\eta),[w,s])| dH_{2n+1}[w,s] dH_{2n+1}(\eta) \\ 
& =  2^\alpha \int_{\partial \cU}  P(\zeta, \eta) \int_{\bH_n} f[w,s] I_\alpha(\Psi(\eta), [w,s])  dH_{2n+1}[w,s] dH_{2n+1}(\eta) \\
& =  2^\alpha P[\widecheck{ I_\alpha f} ](\zeta).
\end{align*}

\end{proof}

\begin{lem} \label{AhernLem1.8}
Let  $ f \in L^1(\bH_n)$ be a positive function. Then there exists a constant $C>0$ such that
\[
M_\gamma (P[\widecheck{I_\alpha f}]) \leq C \widecheck {I_\alpha (f)}.
\]
\end{lem}

\begin{proof}
From Lemma \ref{maximal} we know that $M_\gamma (P[\widecheck {I_\alpha f} ]) \leq C  \widecheck{M (I_\alpha f) } $, so we need to estimate $M (I_\alpha f) $.
\begin{align*}
    MI_\alpha f(x) & = \sup_{r > 0} \frac{1}{r^{2n+2}} \int_{B(x,r)} I_\alpha f(y) dH_{2n+1}(y) \\
    & = C \sup_{r >0} \frac{1}{r^{2n+2}} \int_{B(x,r)} \int_{\bH_n} \frac{f(u)}{d(y,u)^{4(n+1-\alpha)}} dH_{2n+1}(u) dH_{2n+1}(y) \\
    & \leq C \int_{\bH_n} f(u) \sup_{r >0} \big( \frac{1}{r^{2n+2}} \int_{B(x,r)} \frac{dH_{2n+1}(y)}{d(y,u)^{4(n+1-\alpha)}} \big) dH_{2n+1}(u),
\end{align*}
To reach the result it suffices to prove that the supremum is dominated by a constant times $I_\alpha(x,u)$. 

Suppose $d(x,u)^4 \geq 3r$. Since $x \in B(x,r)$, then from triangular inequality it follows that
\[
\frac 1C d(y,u)^4 \leq d(x,u)^4 \leq C d(y,u)^4.
\]
The result follows immediately.

Conversely, if $d(x,u)^4 \leq 3r$, then there exists $C>0$ such that 
\[
B(x,r) \subseteq B(u,Cr).
\]
This yields to
\begin{align*}
    \frac{1}{r^{2n+2}}\int_{B(x,r)} \frac{dH_{2n+1}(y)}{d(y,u)^{2(n+1-\alpha)}} & \leq C \frac{1}{r^{2n+2}} \int_{B(u,Cr)} \frac{dH_{2n+1}(y)}{d(y,u)^{2(n+1-\alpha)}} \\
    & \leq C \frac{1}{r^{2n+2}} \int_0^{Cr} \int_{\partial B(0,1)} \frac{\rho^{2n+1}}{\rho^{2(n+1-\alpha)}} dH_{2n} d\rho \\
    & \leq C \frac{1}{r^{2n+2}} \int_0^{Cr} \rho^{-1+2\alpha} d\rho \\
    & \leq  C \frac{1}{r^{2(n+1-\alpha)}} \\
    & \leq C I_\alpha(x,u).
\end{align*}
\end{proof}

We remark that what we have actually proved in Lemma \ref{AhernLem1.8} is that there exists $C>0$ such that, for every Heisenberg ball $B(x,r)$
\[
\sup_{r>0 }\frac{1}{r^{2n+2}} \int_{B(x,r)} I_\alpha(u,y)dH_{2n+1}(y) \leq C I_\alpha(x,u) \quad \forall x, u \in \bH_n.
\]
In other words the Riesz kernel  $I_\alpha(\cdot, u)$ is a Muckenhoupt weight $A_1$.

\begin{proof}[Proof of the sufficiency of the $\alpha$-subcapacitary condition in Theorem \ref{main_thm}]
    
Let $\mu$ be a subcapacitary measure and consider $F \in H^2_\alpha$. By Lemma \ref{AhernLem1.7}, we know that there exists $f \in L^2(\bH_n)$ such that
\[
|F(\zeta)| \leq P[\widecheck{I_\alpha f}](\zeta),
\]
with $\|F\|_{H^2_\alpha}=\|f\|_{L^2(\bH_n)}$. 

If we use the hypothesis, Theorem \ref{strong_capacitary_inequality} and Lemmas \ref{AhernLem1.7}, \ref{AhernLem1.8} we deduce that 
\begin{align*}
    \int_{\cU} |F(\zeta)|^2 d\mu(\zeta) & = 2\int_0^\infty \mu(\{ \zeta \in \cU: |F(\zeta)| > \lambda\})\lambda d\lambda \\
    & \leq C \int_0^\infty \mu\big(\{ \zeta \in \cU : P [\widecheck {I_\alpha f} ](\zeta) > \lambda \}\big) \lambda d\lambda \\
      & \leq C \int_0^\infty \mu \big( T \{ x\in \bH_n  : \widehat{M_\gamma} (P [\widecheck{I_\alpha f} ])(x) > \lambda \}\big) \lambda d\lambda \\
    & \leq C \int_0^\infty \capa_{\alpha}\big(\{ x \in \bH_n : \widehat{M_\gamma} P [\widecheck{I_\alpha f} ] (x)> \lambda \}\big) \lambda d\lambda \\
    & \leq C \int_0^\infty \capa_{\alpha}\big(\{ x \in \bH_n: I_\alpha f(x) > \lambda \}\big) \lambda d\lambda \\
    & \leq C \|f\|_{L^2(\bH_n)}^2 = C \|F\|_{H^2_\alpha}^2.
\end{align*}
\end{proof}

It is worth noticing that in the proof of the sufficiency part, the fact that $\alpha > \frac{n}{2}$ has not been used in any way. Therefore, the capacitary condition remains sufficient for all $0 < \alpha < \frac{n+1}{2} $. 

We now turn to the necessary part of our main theorem. We shall need a direct estimate of the convolution of the Riesz potential with itself. 

\begin{lem} \label{lem:conv_riesz_kernel}
Let $\frac{n}{2}<\alpha<\frac{n+1}{2}$. Then for every $x,u \in \bH_n$ we have
\[
\int_{\bH_n} I_\alpha(x,y)I_\alpha(y,u)dH_{2n+1}(y) \simeq I_{2\alpha}(x,u).
\]
\end{lem}

\begin{proof}
Let $x, u \in \bH_n$ and set $r=\frac{d(x,u)}{2}$. Then we want to estimate the integral
\[
\int_{\bH_n} \frac{dH_{2n+1}(y)}{d(x,y)^{2(n+1-\alpha)}d(y,u)^{2(n+1-\alpha)}}.
\]
To do that we can split the domain of integration as follows
\begin{align*}
\bH_n & = B(x,r) \cup \{d(x,y) \leq d(y,u)\} \setminus B(x,r) \\ & \cup B(u,r) \cup \{d(y,u) \leq d(x,y) \} \setminus B(u,r).
\end{align*}
We denote by I, II, I', II' the corresponding integrals. By the symmetry of the problem, it is sufficient to estimate I and II.
For the first integral we note that if $y \in B(x,r)$, then by triangular inequality
\begin{align*}
d(y,u) & \leq d(x,u) + d(y,x) \leq 3r \\
d(y,u) & \geq d(x,u) - d(x,y) \geq r.
\end{align*}
So we obtain that
\begin{align*}
    \text{I} & \simeq \frac{1}{r^{2(n+1-\alpha)}}\int_{B(x,r)} \frac{dH_{2n+1}(y)}{d(x,y)^{2(n+1-\alpha)}} \\
    & = \frac{1}{r^{2(n+1-\alpha)}} \int_0^\infty H_{2n+1}\big(\big\{ y \in B(x,r) : d(x,y) \leq t^{\frac{-1}{2(n+1-\alpha)}} \big\}\big) dt.
\end{align*}
We have that $t \leq r^{-2(n+1-\alpha)}$ if and only if $t^{\frac{-1}{2(n+1-\alpha)}} \geq r$. Combining this with \eqref{measure:hball}, we obtain
\begin{align*}
    \text{I} & \simeq \frac{r^{2n+2}}{r^{4(n+1-\alpha)}} \\
    & +  \frac{1}{r^{2(n+1-\alpha)}} \int_{r^{-2(n+1-\alpha)}}^\infty  H_{2n+1}\big(\big\{ y \in B(x,r) : d(x,y) \leq t^{\frac{-1}{2(n+1-\alpha)}} \big\}\big) dt \\
    & \simeq r^{-2(n+1-\alpha)} + \frac{1}{r^{2(n+1-\alpha)}} \int_{r^{-2(n+1-\alpha)}}^\infty t^{-\frac{n+1}{n+1-\alpha}} dt \\
    & \simeq r^{-2(n+1-2\alpha)}.
\end{align*}
We now turn to estimate II.
\begin{align*}
    \text{II} & \leq \int_{\bH_n \setminus B(x,r)} \frac{dH_{2n+1}(y)}{d(x,y)^{4(n+1-\alpha)}} \\
    & = \int_0^\infty H_{2n+1}\big(\big\{ y \in \bH_n : r < d(x,y) \leq t^{\frac{-1}{4(n+1-\alpha)}} \big\}\big) dt.
\end{align*}
Note that if $t > r^{-4(n+1-\alpha)}$ then $t^{\frac{-1}{4(n+1-\alpha)}} < r$, which means
\[
\big\{ y \in \bH_n : r < d(x,y) \leq t^{\frac{-1}{4(n+1-\alpha)}} \big\} = \emptyset.
\]
Since $\frac{n}{2}<\alpha<\frac{n+1}{2}$, then $\frac{n+1}{2(n+1-\alpha)} < 1$ which yields to
\begin{align*}
    \text{II} & \leq \int_0^{r^{-4(n+1-\alpha)}} H_{2n+1}\big(\big\{ y \in \bH_n : r < d(x,y) \leq t^{\frac{-1}{4(n+1-\alpha)}} \big\}\big) dt \\
    & \leq \int_0^{r^{-4(n+1-\alpha)}} H_{2n+1}\big(\big\{ y \in \bH_n : d(x,y) \leq t^{\frac{-1}{4(n+1-\alpha)}} \big\}\big) dt \\
    & \leq \int_0^{r^{-4(n+1-\alpha)}} t^{-\frac{n+1}{2(n+1-\alpha)}} dt \\
    & \simeq r^{-2(n+1-2\alpha)}.
\end{align*}
Combining the two estimates, we obtain that
\[
r^{-2(n+1-2\alpha)} \simeq \text{I} \leq \text{I}+\text{II}  \lesssim r^{-2(n+1-2\alpha)}
\]
which is the desired conclusion.
\end{proof}

In the next lemma we introduce the so called ``holomorphic potentials'' which are holomorphic substitutes of the classical Riesz potentials.

\begin{lem} \label{lem:norm}
    Let $A \subseteq \bH_n$ be a compact set and $\mu^A$ as defined in Theorem \ref{AdamsThm2.5.3}. Consider the {\it holomorphic potential}
    \[
    F_{\mu^A} (\zeta) = \int_{\partial \cU} K_\alpha(\zeta,\omega) d\mu^A(\Psi(\omega)).
    \]
    Then $\|F_{\mu^A}\|^2_{H^2_\alpha} \leq C\capa_\alpha(A)$.
\end{lem}

\begin{proof}
For $r > 0$ and $F \in H^2_\alpha$, we can define the bounded linear functional on $H^2_\alpha$
\[
l_r(F) := \int_{\partial\cU} F(\zeta + \mathbf{i}r) d\mu^A(\Psi(\zeta)),
\]
 where $\mathbf{i}=(0,\dots,0,i)$. Now if we consider $F_r(\zeta):=F_{\mu^A} (\zeta + \mathbf{i}r)$, then we have
 \begin{align*}
\langle K_\alpha^\zeta,F_r \rangle_{H^2_\alpha} & = \ol{F_r(\zeta)} = \int_{\partial \cU} K_\alpha(\omega,\zeta + \mathbf{i}r) d\mu^A(\Psi(\omega))\\ & =\int_{\partial \cU} K_\alpha(\omega +\mathbf{i}r,\zeta) d\mu^A(\Psi(\omega)) \\ 
& = l_r(K_\alpha^\zeta).
 \end{align*}
Since finite linear combinations of the kernels $K_\alpha^\zeta$ are dense in $H^2_\alpha$, then for every $F \in H^2_\alpha$ we have $l_r(F)=\langle F, F_r \rangle_{H^2_\alpha}$. This yields to
\begin{align*}
    \|F_r\|^2_{H^2_\alpha} & = \langle F_r, F_r \rangle_{H^2_\alpha} = l_r(F_r) \\
    & =\int_{\partial\cU} F_{\mu^A}(\zeta + 2\mathbf{i}r)d\mu^A(\Psi(\zeta)) \\
    & = \int_{\partial\cU} \Big( \int_{\partial\cU} K_\alpha(\zeta + \mathbf{i}r, \omega + \mathbf{i}r)d\mu^A(\Psi(\omega)) \Big) d\mu^A(\Psi(\zeta)).
\end{align*}
Since $\Re(K_\alpha) \simeq |K_\alpha|$, we obtain that 
\[
\|F_r\|^2_{H^2_\alpha} \simeq \int_{\partial\cU} \Big( \int_{\partial\cU} |K_\alpha(\zeta + \mathbf{i}r, \omega + \mathbf{i}r)|d\mu^K(\Psi(\omega)) \Big) d\mu^K(\Psi(\zeta)).
\]
Furthermore, note that $| K_\alpha(\zeta + \mathbf{i}r, \omega + \mathbf{i}r)| \leq K_\alpha(\zeta,\omega)$ for $\zeta, \omega \in \cU$. By Lemma \ref{lem:conv_riesz_kernel} and using \eqref{eq:kernel_riesz}, we can conclude that
\begin{align*}
\|F_r\|^2_{H^2_\alpha} & \simeq \int_{\partial\cU} \int_{\partial\cU} I_{2\alpha}(\Psi(\zeta), \Psi(\omega)) d\mu^A(\Psi(\omega)) d\mu^A(\Psi(\zeta))\\
& \simeq \int_{\bH_n} I_\alpha(I_\alpha(\mu^A))(x) d\mu^A(x) \\
& \leq \capa_\alpha(A).
\end{align*}
Combining the above inequality with the dominated convergence theorem, we reach the result.
\end{proof}

\begin{proof}[Proof of the necessity of the $\alpha$-subcapacitary condition in Theorem \ref{main_thm}]
Let $\mu$ be a Carleson measure and let $A = \bigcup_{i=1}^d \overline{B(x_i,r_i)}$ be a finite union of disjoint balls in $\bH_n$. Consider a capacitary measure for $A$, $\mu^A$, as defined in Theorem \ref{AdamsThm2.5.3} and $F_{\mu^A}$ as defined in Lemma \ref{lem:norm}. Combining the hypothesis with Lemma \ref{lem:norm} we reach that
\[
\int_{T(A)} |F_{\mu^A}(\zeta)|^2d\mu(\zeta) \leq \int_{\cU} |F_{\mu^A}(\zeta)|^2d\mu(\zeta) \leq C \|F_{\mu^A}\|^2_{H^2_\alpha} \leq C \capa_{\alpha} (A).
\]
Since $0 \leq \Re F_{\mu^A} \leq |F_{\mu^A}|$, we have that
\begin{equation} \label{prove_subcapacitary}
\int_{T(A)} (\Re F_{\mu^A}(\zeta) ) ^2  d\mu(\zeta) \leq \int_{\cU} |F_{\mu^A}(\zeta)|^2d\mu(\zeta) \leq C \capa_{\alpha}(A).
\end{equation}
We next show that $\Re F_{\mu^A} \geq C$ on $T(A)$. Let $[z,t,h] \in \Psi (T(A))$, then  $B([z,t],h^{1/2}) \subseteq A$. Since by Theorem \ref{AdamsThm2.5.3} we know that $I_{2\alpha}(\mu^A )\geq 1$ a.e. on $A$, we obtain
\begin{align*}
    \Re F_{\mu^A} (\zeta) & = \int_{\partial \cU} P(\zeta,\omega) \Re F_{\mu^A}(\omega) dH_{2n+1}(\omega)\\
    & \simeq \int_{\partial \cU} P(\zeta,\omega) I_{2\alpha}(\mu^A)(\Psi (\omega)) dH_{2n+1}(\omega) \\
    & \geq \int_{\Psi^{-1}(A)} P(\zeta,\omega)  dH_{2n+1}(\omega) \\
    & = C \int_{B([z,t],h^{1/2})} \frac{h^{n+1}}{d([z,t,h],[w,s])^{4(n+1)}} dH_{2n+1}[w,s].
\end{align*}
The task is now to estimate the function $d$. It is easily seen that
\begin{align*}
d([z,t,h],[w,s])^4 & =  \big( \frac{1}{4} |z-w|^2 + 4h \big)^2 + \big( t-s+\frac{1}{2} \Im (z \cdot \overline{w}) \big)^2\\
& \simeq 16h^2 + \frac{1}{16}|z-w|^4 + \big(t-s+\frac{1}{2} \Im (z \cdot \overline{w})\big)^2 \\
& = 16h^2 + d([z,t],[w,s])^4 \\
& \leq 17h^2,
\end{align*}
which yields to
\[
    C\int_{B([z,t],ch^{1/2})}\frac{h^{n+1}}{d([z,t,h],[w,s])^{4(n+1)}} dH_{2n+1}[w,s]  \geq \frac{C}{17^{n+1}} \frac{1}{h^{n+1}} \int_{B([z,t],ch^{1/2})} dH_{2n+1}[w,s]
     \geq C'.
\]
Substituting the above inequalities to \eqref{prove_subcapacitary}, we can conclude that 
\[
\mu(T(A)) \leq C \capa_{\alpha}(A),
\]
which is the desired inequality.
\end{proof}

\bibliography{literature}
\end{document}